\documentclass[12pt,a4paper]{article}
\usepackage[T2A]{fontenc}
\usepackage[cp1251]{inputenc}
\usepackage{amsmath,amssymb,amsthm,amsfonts,amscd}

 \theoremstyle{plain}
\newtheorem{theorem}{Theorem}
\newtheorem*{main*}{Основная Теорема}
\newtheorem*{theorem*}{Теорема}
\newtheorem{lemma}[theorem]{Lemma}

\newtheorem{corollary}[theorem]{Corollary}
\newtheorem{definition}[theorem]{Definition}
\newtheorem{remark}[theorem]{Remark}
\newtheorem*{remark*}{Remark}
\newtheorem*{example*}{Example}

\def\Z{{\Bbb Z}}
\def\R{{\Bbb R}}

\def\Z{{\Bbb Z}}
\def\R{{\Bbb R}}

\def\R{{\Bbb R}}

\begin{document}

\title{Combinatorial formula for the $M$ invariant of magnetic lines}
\author{Akhmetiev P.M.}
\maketitle

\begin{abstract}
	To solve MHD problems within the framework of the theory of two-scale mean fields, it is important to study the invariants of magnetic lines. Such invariants are constructed on the basis of invariants of classical links, which must satisfy the asymptotic property.
	We choose the simplest asymptotic invariant $M_3$ of three-component links, which is not expressed in terms of the pairwise linking coefficients of the components. We check the asymptotic property based only on the combinatorial definition of the invariant and do not use the analytic integral. For simple examples, the proven formula is verified by calculation.
		
\end{abstract}

\section{Introduction}

	Many celestial bodies have  magnetic fields, which, as follows from modern scientific knowledge, are generated by the movement of a conducting liquid medium \cite{Z-R-S}, \cite{Mo}, \cite{R}. This process is known as the dynamo process. The dynamo converts the kinetic energy of the movement of a liquid medium into magnetic energy. The spatial scale of the generated magnetic fields is comparable to the size of the celestial body. Such magnetic fields are called large-scale. The dynamo that generates them is called large-scale.

In dynamo theory, it is generally accepted that the properties of a large-scale dynamo are determined by the behavior of the system on a small scale, which is a hierarchy, and the connection with the large-scale field is described by the mean field equation. The small-scale field is mixed at each scale by the fast field of hydrodynamic velocity of the liquid conducting medium. The large-scale field is the result of averaging the small-scale field. The small-scale field is close to an ideal magnetic field, which is transformed as a magnetic field frozen into a liquid conducting medium. With this approach, the structure of a small-scale field can be characterized by the densities of magnetic line invariants, i.e. such functions that, for an ideal magnetic field, are frozen into the medium and stirred by a fast small-scale velocity field.

It is clear that with this approach there is no question of an exact solution of MHD equations and the statistical properties of the magnetic and hydrodynamic fields come to the fore.  
The concept of MHD spectrum is defined, which characterizes the distribution of hydrodynamic and magnetic energy over scale. In the stationary regime of turbulence, it is necessary to specify a general indicator of the dependence of magnetic and hydrodynamic energy on scale.
In this case, the scale is usually reduced to zero, thereby increasing the wave number to infinity. Of course, the main property of the energy exponent is that the improper integral of the first kind with such an exponent must converge. We will focus on the Kolmogorov MHD spectrum, the properties of which will be recalled in a subsequent publication.
The Kolmogorov MHD spectrum is not the only one; an important role is played by the Iroshnikov-Krechnan spectrum, the properties of which have been studied much better.

We see that for a stable MHD turbulence regime, the average invariant of magnetic lines should be a scale-independent quantity. It is generally accepted that the rate of stirring of magnetic fields is so high that the second derivative rate of dissipation of the magnetic line density invariant becomes very large. In this case, the invariant density is compared to the magnetic flux caused by the $\beta$-term in the mean field equation. Thus, an important role is played by the question: “Are there invariants of magnetic lines, the dimension of which determines a uniform distribution of such an invariant with a given value of the magnetic flux index of the MHD spectrum?”

The dimension of an invariant is described by two real numbers $n_1$, $n_2$, which are written in the form $G^{n_1}sm^{n_2}$. The $n_1$ exponent determines how many times the invariant density value changes if the magnetic field is increased $\mu$-fold? When calculating, it turns out that the exponent $n_1$ is a natural number; as the magnetic field increases, the density of the invariant of magnetic lines increases. The exponent $n_2$ determines how many times the value of the invariant density changes if the scale is increased $\lambda$-fold?  Since we are talking about invariants used for MHD turbulence, the exponent $n_2$ should be negative.   The ratio of the indicators is a wave number that coincides with the magnetic flux indicator of the MHD spectrum. With this value of the magnetic field distribution index, the density of the invariant magnetic lines does not depend on the scale.

It is well known that the linking invariant of magnetic lines satisfies the listed conditions. The coefficient of pairwise linking of magnetic lines in our formulation is characterized by a function on the Cartesian square $\Omega \times \Omega$ of the region with a magnetic field, which is a function of density. Such an invariant exists and is called the asymptotic ergodic Hopf invariant $\chi$, this invariant was discovered by Arnold in 1974 \cite{A-Kh}. The word “asymptotic” means that the invariant is applied to magnetic fields of a general form, the magnetic lines of which are not necessarily closed, but can densely fill part of the volume of the region everywhere, or densely distributed everywhere on the surface.
The word "ergodic" means that the invariant has a density function (since its value is determined by the ergodic theorem applied to magnetic flux) that remains constant during the transformation of the magnetic field region due to hydrodynamics.

The dimension of the asymptotic ergodic invariant is $G^2cm^{-2}$, which means that for a magnetic field distributed with a magnetic flux index of $k=-1$ the density function is independent of scale. It follows from this that the magnetic flux for the Iroshnikov-Kraichnan spectrum with the energy index $k = -\frac{3}{2}$ should, when viewed geometrically, consist of bundles of pairs of magnetic lines that twist among themselves with a constant coupling coefficient $\chi $, independent of scale, normalized to the magnetic length.

Are there invariants of magnetic lines whose dimensions fall under the MHD turbulence indices, for example, under the Kolmogorov energy exponent $k=-\frac{5}{3}$? We give an affirmative answer. Thus, we argue that the magnetic flux of the Kolmogorov spectrum also consists of bundles of quintuplets of magnetic lines that twist and intertwine with each other with a constant value of the M-invariant. This process is already determined not only by the coefficient of twisting of nearby magnetic flux lines, but also by the coefficient of knottedness (intertwining) of nearby magnetic flux lines. The hypothesis about the geometric meaning of the $M$-invariant is formulated.

The invariant $\chi$ (density), the integral convolution of which is called the magnetic helicity invariant, was defined by Gauss on the basis of the Biot-Savart integral. Unlike the $\chi$ invariant, the $M_5$ invariant has not yet been well studied.    It is constructed on the basis of the auxiliary invariant $M_3$, which was discovered in the theory of classical links (a special section of 3D topology) by S.A. Melikhov \cite{M} based on the Conway polynomial in one variable. This invariant has dimension $G^{12}cm^{-6}$ and, thus, has a spectral exponent $k= -\frac{1}{2}$. Such an invariant in magnetic turbulence can only be used for superstrong magnetic fields (the magnetic energy spectrum is distributed in a harmonic series along the wave number scale). If the magnetic field consists only of closed lines, then the $M_3$ invariant is expressed combinatorially as the sum of values for all possible different triplets of magnetic lines.

This is a finite order invariant in the Vassiliev sense. The density value of the $M_3$-integral is given by integral convolution over a finite-multiple configuration of points on all possible triplets of magnetic lines.       This density is calculated in the limit by applying Birkoff's ergodic theorem to the magnetic flux, which generates auxiliary magnetic fluxes in the configuration space.It turns out (this is a new result, the proof of which is not given in the work) that the result of averaging is correct and determines the magnetic line density function, which is defined as an absolutely summable function, almost everywhere on $\Omega^3$, independent of the non-compressive transformation of the domain $\Omega$.

The $M_5$ invariant required to describe the magnetic flux of the Kolmogorov MHD spectrum is determined by the corresponding integral over the configuration space constructed from quintuples of magnetic lines \cite{A1},\cite{A5}. When determining $M_5$ we have to assume that the five magnetic lines are equipped with a fixed cyclic order. One can visually consider that the continuum of magnetic lines is colored with points on the standard circle $S^1$, and the colors are uniformly fractalized in scale.  This kind of hypothesis is a mathematical abstraction that allows one to correctly determine the asymptotic ergodic integral.

From a physical point of view, the presence of $S^1$-colors of magnetic lines is not confirmed.  However, it is possible to determine the symmetrization and for each configuration of five magnetic lines with a prescribed cyclic order, determine the value from all possible recolorings of the magnetic lines in a given configuration. Let us call such symmetrization “loss of symmetry.”  As a result of the loss of symmetry, the invariant $M_5$ is simplified; it turns out that its value can be completely calculated from the values of $M_3$ from all possible triplets of magnetic lines in the chosen quintuplet and from all possible pairwise linking coefficients (the choice of both triplets and pairs in their own subconfigurations can be made unordered).

This work is devoted to the mathematical substantiation of the asymptotic property of the $M_3$ integral.
Research for the $M_5$ integral, most interesting for applications, 
planned for a subsequent publication. In this article we will not be able to answer the question whether the asymptotic property of the link invariant manifests itself 
on an infinite sequence of invariants of increasing complexity, or this property characterizes exclusively the invariants $M_3$, $M_5$.

On the one hand, we are talking about the properties of analytic links, the components of which are linked with positive coefficients and twist to the right. The asymptotic property of arbitrary links was studied in \cite{Ci},\cite{T}, which probably leads to new formulas in applications. On the other hand, invariants with asymptotic properties may turn out to be very rare and may be limited to algebraic constructions with $M$-invariants. The formulated question is interesting not only from a theoretical point of view for the theory of gears. For example, investigations of hierarchical structure
MHD invariants would lead to a new point of view on L.D. Landau’s objection to A.N. Kolmogorov’s statement about the role of the law $\frac{2}{3}$ in hydrodynamics \cite{F}.

Another interesting question is whether the quadratic helicity invariant discussed in \cite{A4},\cite{A-C-S},\cite{E-P-T},\cite{K-V},\cite{D},\cite{C-B} \cite{S}, with higher $M$-invariants in the context of a specific physical problem, for example, in the problems \cite{D-A} related with axion fields?

\section{Combinatorial formula of $M_3$ invariant}

Suppose that $m=3$, consider the oriented link
$\L = L_1 \cup L_2 \cup L_3$. The Melikhov invariant $\gamma(\L)$ is determined by the formula (cf. \cite{A} Theorem 17):

$$
\gamma(\L) = c_1(\L) 
$$
$$-((1,2)(2,3) + (2,3)(3,1) +
(3,1)(1,2))(c_1(L_1) + c_1(L_2) + c_1(L_3))$$
$$ -((3,1)+(2,3))(c_1(L_1 \cup L_2) - (1,2)(c_1(L_1) + c_1(L_2))) $$
$$ -((1,2)+(3,1))(c_1(L_2 \cup L_3) - (2,3)(c_1(L_2) + c_1(L_3))) $$
$$ -((2,3)+(1,2))(c_1(L_3 \cup L_1) - (3,1)(c_1(L_3) + c_1(L_1))), $$
where  $(i,j)=c_0(L-i \cup L_j)$ is the linking number of the pair of components  $L_i$, $L_j$, $i,j = 1,2,3$, $i \ne j$,
$c_1$ is the coefficient in the Conway polynomial 
 $\nabla_{\L}(z)=z^{m-1}(c_0 + c_1z^2 + \dots )$
for the corresponding proper $m$--component link. The work uses standard equipment 
calculations of the Conway polynomial, it is presented, for example, in \cite{P-S}; Our calculations are close to those from \cite{M}.

.

Let us define the invariant $\tilde M(\L)$ by the formula:
\begin{eqnarray}\label{MInv}
	\begin{array}{c}
\tilde M(\L) = -(1,2)(2,3)(3,1)\gamma(\L) \\
+ ((1,2)^2(1,3)^2\beta(L_2 \cup L_3) + (2,3)^2(2,1)^2\beta(L_3 \cup
L_1) + (2,3)^2(2,1)^2\beta(L_3 \cup L_1)),\\
\end{array}
\end{eqnarray}
where
\begin{eqnarray}\label{SatoLevine}
	\beta(L_i \cup L_j) = c_1(\L) - c_0(\L)(c_1(L_1) + c_1(L_2))
\end{eqnarray}
is the Sato-Levine invariant for the oriented 2-component link
$\L = L_i \cup L_j, \quad i \ne j$.


\begin{theorem}\label{th1}
	1. There is a symmetric odd polynomial $P((2,3),(3,1),(1,2))$ (see the formula $(\ref{Poly3})$ ) such that
	\begin{eqnarray}\label{m}
		M(\L) = \tilde{M}(\L) + P((2,3),(3,1),(1,2))
	\end{eqnarray}
	is asymptotic. This means that for an arbitrary 3-vector
	$(3,1),(1,2),(2,3)$ with non-negative (non-positive, respectively) coordinates and for an arbitrary vector $(\lambda_1,\lambda_2,\lambda_3)$
	with positive coordinates
 the invariant $M$ satisfies the relation: 
	\begin{eqnarray}\label{ass}
	M(\lambda_1 L_1 \cup \lambda_2 L_2 \cup \lambda_3 L_3) = \lambda_1^4 \lambda_2^4 \lambda_3^4 M(L_1 \cup L_2 \cup L_3).
	\end{eqnarray}
	(more detailed definition   $L_1 \cup L_2 \cup L_3 \mapsto \lambda_1 L_1 \cup \lambda_2 L_2 \cup \lambda_3 L_3$,
	$\lambda_i \in \Z_+$, 
	$i=1,2,3$ will be in Section  4).

	2. For the standard (left) Hopf link  $L_{Hopf}$ с $(1,2)=(2,3)=(3,1)=-1$
	(definition in Section  \ref{Sec3}), the following equation is satisfied: 
	\begin{eqnarray}\label{n}
		 \tilde{M}(\L_{Hopf}(-1,-1,-1)) = +1.
	\end{eqnarray}


	3. The combinatorial invariant in the right-hand side of the foremula  $(\ref{m})$
	coinsids with the integral invariant from \cite{A}. 
	
	4. The invariant $M$ in (\ref{m}) is a mirror-symmetric. This means that the mirror symmetry (the mirror symmatry changes orientations on the components of $\L$ and is a reflection with respect to a plane in the space; the mirror symmetry changes the collection of pairwise linking numbers of a link into the opposite):
	$\{(2,3),(3,1),(1,2)\} \mapsto \{-(2,3),-(3,1),-(1,2)\}$)
	changes $M \mapsto -M$ into the opposite.

	5. The invariant $M$  is extended for links with an arbitrary collection of pairwise linking number, the extended invariant is not changed with respect to the following inversion of orientations on a pair of components:  $L_j \mapsto -L_j; L_{j+1} \mapsto -L_{j+1}, L_{j+2} \mapsto L_{j+2}$.

	6. The polynomial in $(\ref{m})$ (for a right link with positive linking numbers) is defined by the following formula:
	\begin{eqnarray}\label{Poly3}
		P=P_1+R 		
		\end{eqnarray}
	where $P_1$ is defined by (\ref{ppoly1}), $R$ is defined by the formula (\ref{R}).
\end{theorem}

\begin{remark*}

The integral invariant $M$, constructed in \cite{A} does not depend on an order of magnetic lines and of  orientations of magnetic lines in the triple (see statement 3 of Theorem \ref{th1}). By the construction an inversion of an orientation on a one component of the oriented link corresponds to permutation of the order of the line by the corresponding odd transposition (see the subsection "Internal symmetry"). We will not investigate the gauge of such a transformation. In the prove that $M$ is a Vassiliev finite-type invariant we use
the configuration space of the ordered oriented magnetic lines. (In particular, the formula (19),\cite{A}
changes into the opposite, when we apply an odd transposition of an order of the components).

\end{remark*}

\begin{remark*}
	The standard (left) Hopf link consists of three layers of the Hopf fibration, with
	the orientation in the bundle space is chosen in such a way that all pairwise linking coefficients of components are equal to $-1$ (see Fig. 1; cf. Fig. for a right Hopf link).

	\[   \]
	Fig. 1
	\[   \]
	
	Equivalently, the bundle on the standard sphere $S^3$ is chosen to be the Hermitian conjugate to the holomorphic Hopf bundle. The space $S^3 \setminus \{pt\}$ is identified with the $\R^3$ stereographic projection, which allows us to write: $\L_{Hopf}(-1,-1-1) \subset \R^3$ .

\end{remark*}

\begin{remark*}
	In \cite{A} a conjecture is made about the polynomial $P$ in the formula (\ref{m}). The conjecture is proven in paragraph 3 of Theorem \ref{th1}, where it is stated that the polynomial $P$ coincides with the polynomial $P_1$ defined below by the formula (\ref{ppoly1}). We refute this hypothesis.
	\end{remark*}

\begin{remark}
	Let's apply the formula \ref{m} for a new way to calculate the Casson invariant of a knot $K \subset \R^3$ \cite{M-P}. 
	Consider a 3-component satellite link $\L$ associated with $K$ whose components are linked to each other with a nonzero 
	coefficient $p$ ($p=-2$ on the figure).

	\[  \]
	Fig. 2
	\[  \]
	
	Note that the terms $c_1(K)$ included in the formula (\ref{MInv}) is not canceled.
	For $p=-2$, for example, it turns out:
	$\tilde{M}(K)=680$. For the link $\L_{Hopf}(+2,+2,+2)$ in the mirror $\Delta$-equivalent class of $\L$ we get $\tilde{M}=-104$ (comp. with calculations attached to Fig.5).
\end{remark}

An attempt to construct an asymptotic ergodic integral formula for magnetic knots based on the Sato-Lewin $\beta$ invariant was made in \cite{A-K}. It turned out that the formula works only in the case when the coefficient of engagement of two components of engagement is equal to zero. In this case, the $\beta$ invariant applied to the two-component satellite link of a node is not related to the knottedness of this knot.

The invariant formula $M$ for a three-component link works for an arbitrary three-component link.  It is impossible to apply such a formula without additional constructions to calculate the invariants of magnetic knots, since the formula has a peculiarity when the components of the satellite linkage approach each other.

The formula we obtained using dipole sources for the vector potential provides a new integral formula for the Casson invariant of a knot. Such an integral is absolutely divergent, but must converge in the sense of the principal value. Informally speaking, somewhere in the depths of the configuration space there are multipoles, which contribute to the main value of the integral for the Casson invariant.

\section{Special 2- and 3-component oriented links }\label{Sec3}

We need to define special 2-component links $\L_{Hopf}(p)$, $\L^{op}_{Hopf}(p)$ (with the order of the components)
each of which depends on the integer parameter $p \in \Z$ - the coefficient of linkage of the components.
Special 3-component links $\L_{0}((2,3),(3,1),(1,2))$ (with the order of the components) will be required, see Fig. 3.

\[  \]
Fig.3
\[  \]

The corresponding pair of components of this link is $\L_{Hopf}^{op}((3,1))$ as in Fig. 4.
\[  \] 
Fig.4
\[  \]
Let us also define the 3-component link $\L_{Hopf}(p,p,p)$, see Fig. 5,
\[  \]
Fig.5
\[  \]
depending on an integer parameter (this parameter is equal to the pairwise coefficient of engagement of the components, which are all the same).  Let us define modified links
\begin{eqnarray}\label{L1}
\begin{array}{c}
\L_{0}((2,3),(3,1),(1,2)) \mapsto \L_{0}^{op(1)}((2,3),-(3,1),-(1,2)),\\
\L_{0}((2,3),(3,1),(1,2)) \mapsto \L_{0}^{op(2)}(-(2,3),(3,1),-(1,2)),\\ 
\L_{0}((2,3),(3,1),(1,2)) \mapsto \L_{0}^{op(3)}(-(2,3),-(3,1),(1,2)).
\end{array}
\end{eqnarray}
\begin{eqnarray}\label{L2}
\L_{Hopf}(p,p,p) \mapsto \L_{Hopf}^{op(1)}(p,-p,-p).
\end{eqnarray}

Let's start with the definition
$\L_{Hopf}^{op}(p)$. First component $L_1$ of the link $\L_{Hopf}^{op}(p)$
is a standard flat circle, the second component is $L_2$
is located on the boundary of the tubular neighborhood, this component rotates once along the meridian along the central line, but oriented in the opposite direction with respect to the first component; around the parallel
the component rotates $p$ times and the linkage coefficient of the components is equal to the prescribed integer value $p$.


Let's define a link 
$\L_{Hopf}(p)$, which is obtained from
$\L^{op}_{Hopf}(-p)$ by reversing the orientation of the second component. Note that when the orientation of both components is simultaneously reversed, the isotopic class
$\L_{Hopf}(p)$.
Component engagement coefficient
We will denote $\L_{Hopf}(p)$ by $p$. From an algebraic point of view, the link 
$\L_{Hopf}^{op}(-p)$ is simpler compared to $\L_{Hopf}(p)$,
since $c_1(\L_{Hopf}^{op}(-p))=0$, see Lemma \ref{lemma10}.

Let us define $\L_{Hopf}(p,p,p)$ as the link obtained by resolving a plane circle with self-linking coefficient $p$ into three link components pairwise linked with coefficient $p$. For $p=1$ the definition is required in Theorem \ref{th1} item 2. Modifications (\ref{L1}), (\ref{L2}) are defined in an obvious way by reversing the orientation of the component marked with a superscript in the notation.

\section{Calculation of Conway invariants for special links}

The following lemma on the coefficient 
 $c_1(\L)$ in the Conway polynomial
for two-component links $\L = L_1 \cup L_2$. In the case components  $L_1$, $L_2$ are unknotted, by (\ref{SatoLevine}) we get:
$c_1(\L)=\beta(\L)$.

\begin{lemma}\label{lemma10}
	
	Sato-Levine invariant, determined by the formula
	$(\ref{SatoLevine})$, is calculated by the following formula: 
	
	\begin{eqnarray}\label{hopf}
		\beta(\L_{Hopf}(p))= \frac{(p+1)p(p-1)}{6}.
	\end{eqnarray}
\end{lemma}

\subsection*{Proof of Lemma $\ref{lemma10}$}

Let us offer two different proofs.

First proof.
Calculating the left side of the formula
$(\ref{hopf})$ using the Conway polynomial proves that 
$\beta(\L_{Hopf}(p))$ is a third-degree polynomial in the variable $p$.
For $p=-1,0,$ or $+1$, we get 
$\beta(\L_{Hopf}(p))=0$ and $\beta(\L_{Hopf}(2))=1$.
The equality $(\ref{hopf})$ is the only possible one.

Second proof.
To engage
$\L_{Hopf}(p)$ consider the link
$\L^{op}_{Hopf}(-p)$. As we noted above, it will turn out
$c_1(\L^{op}_{Hopf}(-p))=0$. The components of $\L^{op}_{Hopf}(-p)$ are unknotted
and according to the formula
$(\ref{SatoLevine})$ we get $\beta(\L^{op}_{Hopf}(-p))=0$.
Formula $(\ref{hopf})$
is a special case of the formula for the case when one of the components of the link changes
orientation. In this general formulation, the formula is proven in
\cite{Ni}. \qed
\[  \]

\begin{lemma}\label{Hopf3}
	The following equalities are valid:
	$$ \gamma(\L_{0}((2,3),(3,1),(1,2)))=0, $$
	$$ \tilde{M}(\L_0((2,3),(3,1),(1,2))=0. $$
\end{lemma}

\subsubsection*{Proof of the Lemma \ref{Hopf3}}
Let's prove the first equality. The components of the link $\L_{0}((2,3),(3,1),(1,2))$ are unknotted and the proper two-component sublinks are links $\L_{Hopf}^{op}(k)$, 
for $p=(2,3), p=(3,1), p=(1,2)$, for which 
$c_1(\L_{Hopf}^{op}(p))=0$.
Consequently, the following equality holds for the Melikhov invariant:
$$ \gamma(\L_{0}((2,3),(3,1),(1,2)))=c_1(\L_{0}((2,3),(3,1) ,(1,2))). $$
Let's begin to untangle the link by reducing the absolute value of the coefficient $(1,2)$ to zero, using the homotopy with the intersection of $L_1$, $L_2$.
The next step of the design will result in a gear
of the same type, for which the coefficient $(1,2)$ is $1$ less in absolute value.
We will prove the required equality by induction.

In this case, we will use the relation for the Conway polynomial.
At each intersection of the $L_1$ and $L_2$ components of the singular link, which we denote by  
$L_1 \sharp L_2 \cup L_3$, the number of components decreases 
from $3$ to $2$, so the jump in the value of $c_1(\L)$ is determined by the value of $c_1(L_1 \sharp L_2 \cup L_3)$.
But $L_1 \sharp L_2 \cup L_3 = \L_{Hopf}^{op}((1,3)+(2,3))$. Since $c_1(\L_{Hopf}^{op}(p))=0$, $p \in \Z$,
we obtain the required induction step.

The basу of еру induction is to check equality:
$$c_1(\L_{0}((2,3),(3,1),0)) = 0.$$ 
This equality can be proved in a similar way. Used 
additional induction that reduces the absolute value of the coefficient $(2,3)$ to zero.
Each step of additional induction is similar to the step of the main one; the basis of additional induction follows from the fact that for a three-component link obtained from $\L^{op}_{Hopf}((3,1))$, adding a disjoint
of the unknotted component $L_3$, the value $c_1$ becomes zero. Finally, the basis of the second induction follows from the relation $c_1(\L_0(0,(3,1),0))=0$, which follows from the fact that the link $\L_0(0,(3,1),0) $ has an unknotted component $L_2$ in a ball that does not intersect components $L_1, L_3$.

The second equality easily follows from the formulas $\beta(\L_{Hopf}^{op}(p))=0$, $c_1(\L_{Hopf}^{op}(p))=0$.
\qed

Consider the 3-component link $\L_{Hopf}(-1,-1,-1)$, which is used in the formulation of Theorem \ref{th1}.

\begin{lemma}\label{Hopf4}
	The following equation is satisfied:
	$$\gamma(\L_{Hopf}(-1,-1,-1))=+1. $$
\end{lemma}

\subsubsection*{Proof of the Lemma \ref{Hopf4}}

As in the proof of Lemma \ref{Hopf3}, we will untangle the link $\L_{Hopf}(-1,-1,-1)$
homotopy of components $L_2$, $L_3$ with a self-intersection point, ensuring that these components become unlinked. Let us denote disentanglement by
$$\L_{Hopf}(-1,-1,-1) \mapsto \L_{Hopf}(0,-1,-1). $$
The equality is satisfied: 
$c_1(\L_{Hopf}(0,-1,-1))=0$, which was proven in Lemma \ref{Hopf3}.

At a single intersection of the components $L_1$ and $L_2$, a singular link arises, which we denote by  
$L_1 \sharp L_2 \cup L_3$. The jump in the value of $c_1(\L)$ is determined by the value of $-c_1(L_1 \sharp L_2 \cup L_3)$ since the sign of the intersection is negative.

But $L_1 \sharp L_2 \cup L_3 = \L_{Hopf}(-2)$. Since according to Lemma \ref{lemma10} the equality is true: $c_1(\L_{Hopf}(-2))=-1$,
we get what we need. \qed

Consider the link $\L_0((2,3),(3,1),(1,2))$ and modify this link
$$ \L_0((2,3),(3,1),(1,2)) \mapsto \L_0^{op(1)}((2,3),-(3,1),-( 1,2)). $$
For the integral invariant $M$, by \cite{A} 
the identity is true:
\begin{eqnarray}\label{Mop}
	M(\L_0((2,3),(3,1),(1,2))) = M(\L_0^{op(1)}((2,3),-(3,1),-(1,2)). 
\end{eqnarray}
Let us calculate the value $\tilde{M}$ of the modified link.

\begin{lemma}\label{Hopf6}
	The following equation is satisfied:
	$$ \tilde{M}(\L_0^{op(1)}((2,3),-(3,1),-(1,2)) = $$
	$$-\frac{1}{3}(1,2)^2(2,3)^2(3,1)^2[(1,2)+(3,1)] $$
	$$ -\frac{1}{6}(1,2)(2,3)(3,1)[(1,2)(2,3)+(2,3)(3,1)].$$
\end{lemma}

\begin{corollary}\label{1cor}
	For an arbitrary link $\L$ with pairwise link coefficients $(1,2),(2,3),(3,1)$
	the correct formula is:
	\begin{eqnarray}\label{for1}
		\begin{array}{c}
		\tilde{M}(\L^{op(1)})-\tilde{M}(\L)= \\
		-\frac{1}{3}(1,2)^2(2,3)^2(3,1)^2[(1,2)+(3,1)] 
		 -\frac{1}{6}(1,2)(2,3)(3,1)[(1,2)(2,3)+(2,3)(3,1)].
		 \end{array}
		\end{eqnarray}
\end{corollary}

\subsubsection*{Proof of Lemma  \ref{Hopf6}}
Let us prove the equality:
\begin{eqnarray}\label{c1}
	\begin{array}{cc}
		6c_1(\L_0^{op(1)}((2,3),-(3,1),-(1,2))) = -(2,3)[((1,
		2)+(3,1))^3-((1,2)+(3,1))] \\
		+ (1,2)[(1,3)^3-(1,3)] + (3,1)[(1,2)^3-(1,2)].
	\end{array}
\end{eqnarray}
Let us begin to untangle the components $L_2, L_3$ by a homotopy with a self-intersection point, ensuring that these components are unlinked. At the next step, if $(2,3)\ne 0$, we obtain a link of the same type, whose coefficient $(2,3)$ is $1$ less in absolute value.
The base of the induction follows from the following equality:

\begin{eqnarray}\label{op}
	6c_1(\L_0^{op(1)}(0,-(3,1),-(1,2)))= (1,2)[(1,3)^3-(1,3)] + (3,1)[(1,2)^3-(1,2)].
\end{eqnarray}

The equality $(\ref{op})$ is proved by the same methods as Lemma \ref{Hopf4}.
Consider the homotopy of the link $\L_0^{op}(0,-(3,1),-(1,2))$, decoupling the components $L_1$ and $L_3$, fixed outside the ball containing the component $L_2$. The homotopy has $(3,1)$ points (taking into account the sign) of self-intersection.
As a result of such a homotopy we obtain a link $\L_0^{op}(0,0,-(1,2))$ for which $c_1(\L_0^{op}(0,0,-(1,2)) )=0$. The value $c_1(\L_0^{op}(0,-(3,1),-(1,2)))$ is calculated as the sum of all values
$2$-component links $\L(i)=L_2 \cup L_3 \sharp L_1(i)$ arising in the process of homotopy.

The equality is valid: $c_1(\L(i))=(3,1)c_1(L_3 \sharp L_1(i)) + c_1(\L_{Hopf}((1,2)))$,
which can be proven by calculating $\beta(\L(i)) = \beta(\L_{Hopf}((1,2))$, then using the formula
$c_1(\L(i))=\beta(\L(i)) - (3,1)c_1(L_3 \sharp L_1(i))$.

Справедливо равенство: $c_1(\L(i))=(3,1)c_1(L_3 \sharp L_1(i)) + c_1(\L_{Hopf}((1,2)))$,
которое можно доказать вычисляя $\beta(\L(i)) = \beta(\L_{Hopf}((1,2))$, затем используя формулу
$c_1(\L(i))=\beta(\L(i)) - (3,1)c_1(L_3 \sharp L_1(i))$.
But
$$ \sum_i c_1(L_3 \sharp L_1(i)) = \beta(\L_{Hopf}(3,1)) = c_1(\L_{Hopf}((3,1)). $$
Now it is enough to use the formula $(\ref{hopf})$, you get the formula $(\ref{op})$ for the base of induction.

At an each step of the induction we get the following link:
$-L_1 \cup L_2 \sharp L_3 = \L_{Hopf}^+(-(1,2)-(3,1))$,
for which by the formula
 (\ref{hopf}) we get:
$6c_1(\L_{Hopf}^-(-(1,2)-(3,1)))=-((1,
2)+(3,1))^3-((1,2)+(3,1))$.
The equality
 $(\ref{c1})$ is proved.

Let us calculate the values of the remaining coefficients of the Conway polynomial included in the formula $\gamma$.
Obviously,
 $c_1(L_1)=c_1(L_2)=c_1(L_3)=0$, $c_1(L_2 \cup L_3)=0$. 
We have: $6c_1(-L_1 \cup L_2)=-((1,2)^3-(1,2))$, $6c_1(-L_3 \cup L_1)=-((3,1)^3-(3,1))$.
Let us calculate the value $6\gamma(\L_0^{op}((2,3),-(3,1),-(1,2)))$. We get:
$$-(2,3)[((1,2)+(3,1))^3-((1,2)+(3,1))]  + (1,2)[(1,3)^3-(1,3)] + (3,1)[(1,2)^3-(1,2)]$$
$$ + (-(3,1)+(2,3))[(1,2)^3-(1,2)]  
+((2,3)-(1,2))[((3,1)^3-(3,1)]. $$

Let us calculate the value
$6\tilde{M}(\L_0^{op}((2,3),-(3,1),-(1,2)))$. We get:
$$-(1,2)(2,3)^2(3,1)[((1,2)+(3,1))^3-((1,2)+(3,1))]] $$
$$+ (1,2)^2(2,3)(3,1)[(1,3)^3-(1,3)] + (1,2)(2,3)(3,1)^2[(1,2)^3-(1,2)]]  $$
$$ +(1,2)(2,3)(3,1)[-(3,1)+(2,3)][(1,2)^3-(1,2)] $$
$$+(1,2)(2,3)(3,1)[(2,3)-(1,2)][((3,1)^3-(3,1)]  $$
$$+ (2,3)^2(2,1)^2[(3,1)^3-(3,1)] + (2,3)^2(2,1)^2[(1,2)^3-(1,2)].$$

This equality is equivalent to the required equality. \qed



\section{The Polynomial (\ref{m}) }\label{Sec4}



The integral invariant $M(\L)$ does not change when the components of the link are renamed and is skew-invariant under mirror symmetry of the space. 
The combinatorial invariant $\tilde{M}$ also does not change when the components of the link are renamed and is skew-invariant under mirror symmetry of the space.

It follows from this that
the polynomial $P$ must be expressed as a sum of homogeneous symmetric odd polynomials in three variables:
\begin{eqnarray}\label{Poly}
	 P((2,3),(3,1),(1,2))=P_1((2,3),(3,1),(1,2))+R((2,3),(3,1),(1,2)),
\end{eqnarray}	
where the polynomial $P_1((2,3),(3,1),(1,2))$ consists of terms (\ref{case1}), and $R((2,3),(3,1) ,(1,2))$ consists of terms (\ref{case2}).

The polynomial (\ref{Poly}) is presented as a sum of symmetric homogeneous elementary monomials:
\begin{eqnarray}\label{poly1}
	P((2,3),(3,1),(1,2)) = \sum_{i=0}^N P_i((2,3),(3,1),(1,2)). 
\end{eqnarray}

The monomials in a symmetric polynomial can have the form $(1,2)^i(2,3)^j(3,1)^k$, where the following two cases are a priori possible (up to a cyclic permutation).

Case 1: 
\begin{eqnarray}\label{case1}
i=1, j=0, k=0 \quad \pmod{2}. 
\end{eqnarray}

Case 2. 
\begin{eqnarray}\label{case2}
i=1,j=1,k=1 \quad \pmod{2}.
\end{eqnarray}

\begin{lemma}\label{l17}
	If at least one coefficient $(2,3), (3,1), (1,2)$ becomes zero,
	$(2,3)(3,1)(1,2)=0$, then the polynomial $P$ vanishes identically.
\end{lemma}

\begin{corollary}
	From Lemma \ref{l17} the first three terms in $(\ref{poly1})$, for $i=0,1,2$, vanish.
\end{corollary}

\subsubsection*{Proof of Lemma \ref{l17}}
If
two of the three coefficients vanish,
$$ (1,2)^2(2,3)^2 + (2,3)^2(3,1)^2 + (3,1)^2(1,2)^2 = 0, $$
then $\tilde{M}(\L)$ in the right-hand side of the formula $(\ref{m})$ is zero. 
But in this case the integral invariant
$M(\L)$ also vanishes, as easily follows from \cite{A1}, \cite{A3}.

If only one gearing coefficient is zero, and the other two are non-zero, you get:
$$ (1,2)^2(2,3)^2 + (2,3)^2(3,1)^2 + (3,1)^2(1,2)^2 \ne 0, \quad (1,2)(2,3)(3,1)=0. $$
Let's say, for definiteness, $(1,2)(2,3) \ne 0$,
then using the formula $(\ref{m})$ we get:
$$\tilde{M}(\L) = (1,2)^2(2,3)^2 \beta(L_1 \cup L_3). $$

But in this case, according to \cite{A3} it will also turn out
$$ M(\L) \simeq (1,2)^2(2,3)^2 \beta(L_1 \cup L_3), $$
Moreover, the normalization coefficient in this equality does not affect the value of the polynomial
$P((2,3),(3,1),(1,2))$, since in this case the polynomial vanishes.
 \qed

\begin{lemma}\label{l19}
	Put the following polynomial (\ref{ppoly1}) in the formula (\ref{m}). 
\begin{eqnarray}\label{ppoly1}
	\begin{array}{c}
		P_1((2,3),(3,1),(1,2))=                                         \\
		+\frac{1}{6}(1,2)^2(2,3)^2(3,1)^2[(1,2)+(2,3)+(3,1)]            \\
		+\frac{1}{12}(1,2)(2,3)(3,1)[(1,2)(2,3)+(2,3)(3,1)+(3,1)(1,2)].
	\end{array}
\end{eqnarray}
Then the formula (\ref{ass}) is satisfyed for the following vectors:
\begin{eqnarray}\label{v1}
(\lambda_1,\lambda_2,\lambda_3)=(+1,-1,-1),(-1,+1,-1),(-1,-1,+1);
\end{eqnarray}
\begin{eqnarray}\label{v2}
(\lambda_1,\lambda_2,\lambda_3)=(-1,+1,+1),(+1,-1,+1),(+1,+1,-1),(-1,-1,-1).
\end{eqnarray}

This exactly means that, $\tilde{M}+P_1$ is an invariant of non-oriented links.
\end{lemma}

\subsubsection*{Proof of Lemma \ref{l19}}
Let $\L=\L((2,3),(3,1)(1,2))$ be a link with prescribed pairwise link coefficients.
In particular, the case is used when $\L=\L_0$ is a special link.
Let's consider modifications
$\L^{op(3)}=\L^{op(1)}(-(2,3),-(3,1),(1,2))$, $\L^{op(1)}=\L^{op(2)}((2,3),-(3,1),-(1,2))$, $\L^{op(2)}=\L^{op(3)}(-(2,3),(3,1),-(1,2))$, $\L^{op(1+2+3)}((1,2),(2,3),(3,1))$. 

Let's determine the average value: 
\begin{eqnarray}\label{av1}
 \tilde{M}^{av}(\L) = \frac{1}{4}[ \tilde{M}(\L^{op(1+2+3)}) + \tilde{M}(\L^{op(3)}) + \tilde{M}(\L^{op(1)}) + \tilde{M}(\L^{op(2)})].
 \end{eqnarray}

Because $\tilde{M}$ is defined using Conway polynomial, which is invariant with respect to inversion of the orientation of all components of the link, we get:
$\tilde{M}(\L^{op(1+2+3)})=\tilde{M}(\L)$.

It is clear that the value of $\tilde{M}^{av}(\L)$ is the same for any modification,
described in the formula (\ref{v1}). Values of invariants in the two orbits of the modifications (\ref{v1})(\ref{v2}) equal correspondingly. Lemma \ref{l19} is proved. \qed

\subsubsection*{Good links}
Suppose that $(1,2)(2,3)(3,1) \neq 0$.
Consider an arbitrary link $\L=L_1 \cup L_2 \cup L_3$ and define the normalization transformation
$$ (L_1 \cup L_2 \cup L_3) \mapsto \L_{norm}, $$
which is the result of cable construction
$$ (L_1 \cup L_2 \cup L_3) \mapsto ((2,3)L_1,(3,1)L_2,(1,2)L_3) = \L_{norm}. $$
With an additional orientation reversal transformation on the first component:
$$ ((2,3)L_1,(3,1)L_2,(1,2)L_3) = \L_{norm} \mapsto ([(2,3)L_1]^{op},(3,1)L_2,(1,2)L_3^{op})=  \L_{norm}^{op}.$$ 
For a link $\L_{norm} = ((2,3)L_1,(3,1)L_2,(1,2)L_3)$ all pairwise link coefficients are identical and equal to the product $k=(1,2)(2 ,3)(3,1)$. For the link $\L_{norm}^{op}=([(2,3)L_1]^{op},(3,1)L_2,(1,2)L_3)$ the pairwise link coefficients are equal:
 $(k,-k,-k)$, where $k=(2,3)(3,1)(1,2)$.

We say that a link $\L$ with nonnegative pairwise link coefficients is good if for its pairwise link coefficients
two properties are satisfied.

1. coefficients are perfect squares.

2. system of equations:
$$ \sqrt{(1,2)} = \mu_1 \mu_2,\quad \sqrt{(2,3)} = \mu_2 \mu_3, \quad \sqrt{(3,1)} = \mu_3 \mu_1 $$	
has a positive integer solution
$(\mu_1,\mu_2,\mu_3)$.

A good link with positive pairwise link coefficients has the same set
pairwise linkage coefficients, as the linkage
$ \L_{Hopf}(\mu^2_1,\mu^2_3,\mu^2_2)$ obtained by laying a secondary $(\mu_1,\mu_2,\mu_3)$-cable
for the primary cable $L_{Hopf}(\mu_1,\mu_2,\mu_3)$ of the Hopf link $\L_{Hopf}(+1,+1,+1)$. If the signs of the gearing coefficients are arbitrary, then the definition of good gearing is formulated as follows.

For links, both right and left, the signs $(1,2), (2,3),(3,1)$ correspond 
that
a corresponding secondary $(\mu_1,\mu_2,\mu_3)$-cable is laid to engage $\L_{Hopf}( \mu_1 (\pm1),\mu_2(\pm1),\mu_3(\pm1))$, for whose signs of pairwise linking coefficients are chosen as the signs of $(2,3),(3,1),(1,2)$.

For a good (right) link $\L$ we denote by $\sqrt {k } = \mu^2_1\mu^2_2\mu^2_3$ the positive root.
In each case, right and left, we define $k = (\sqrt{k})^2 > 0$. 
Note, however, that for a right link $(1,2)(2,3)(3,1)>0$, for a left link $(1,2)(2,3)(3,1)< $0. The $k$ sign was already taken into account when laying the cable as a sign of a special Hopf link.

\subsection*{Calculation for $\L_{Hopf}(k,k,k)$}
Let $p>0$ be a parameter.
Consider the link $p\L_{Hopf}(+1,+1,+1)$ obtained by laying the cable of a standard right Hopf link.

Consider a sequence of $\Delta$-movements connecting 
$p\L_{Hopf}^{op(1)}(+1,-1,-1)$ and $\L_{Hopf}^{op(1)}(p^2,-p^2,- p^2)$, let us study the jump of the $\tilde{M}$-invariant. This sequence of $\Delta$-movements is defined since the values of the pairwise coefficients of engagement of the components of this pair of engagements coincide.

\begin{lemma}\label{l23}
	Jump of the invariant $\tilde{M}$ during a sequence of $\Delta$-motions on a pair of links
	$p\L_{Hopf}^{op(1)}(+1,-1,-1)$ and $\L_{Hopf}^{op(1)}(p^2,-p^2,- p^2)$ is equal to zero.	
\end{lemma}

\begin{lemma}\label{l24}
	The following equalities are valid:
	\begin{eqnarray}\label{eql24}
		\begin{array}{c}
		\tilde{M}(\L_{Hopf}(p,p,p)) = +\frac{p^7}{2}+\frac{p^5}{2 } \\
		\tilde{M}(\L_{Hopf}^{op(1)}(p,-p,-p)) = -\frac{p^{7}}{6} + \frac{p^{5}}{6}.
		\end{array}
	\end{eqnarray}
\end{lemma}

\begin{corollary}
	The equality is satisfied:
	\begin{eqnarray}\label{eql24}
		\begin{array}{c}
			\tilde{M}^{av}(\L_{Hopf}(p,p,p)) = \tilde{M}^{av}(\L_{Hopf}^{op(1)}(p,-p,-p))= \frac{1}{4}p^{5}.
		\end{array}
	\end{eqnarray}
\end{corollary}

\subsubsection*{Proof of Lemma  \ref{l23}}

Consider the link $\L_{Hopf}^{op(1)}(p^2,-p^2,-p^2)$ and bring together the components $L_2,L_3$, which we place inside the thin solid torus $U$ as two $p$-cable of Hopf link
$(L_2 \cup L_3)=\L_{Hopf}(+1) \subset U$. Component $L_1$ is located in the complement of $U$ and the view of this component is
is unimportant, since under the homotopy $L_1 \mapsto L_1'$ inside $\R^3 \setminus U$ of this component the $\tilde{M}$-invariant of the link $(L_1 \cup L_2 \cup L_3)$ does not change.

Consider a circle $\l \subset U$ that wraps $p$-times around the center line
$S^1 \subset U$ of a solid torus and consider the frame $\xi$ of this circle $\l$ with a self-linking coefficient $p^2$, for which the framed knot $(\l,\xi)$ is isotopic to the standard flat circle with the same self-linking coefficient $p^2$.
Let us define $L'_2 \cup L'_3 \subset U$ as a result of shifting $l$ along the framing vector $\xi$. By design the sublink
$L'_2 \cup L'_3$ isotopically is $\L_{Hopf}(p^2)$. Transform by homotopy $L_1 \mapsto L'_1$
component so that $(L_1 \cup L_2 \cup L_3)$ is isotopic to $\L_{Hopf}^{op(1)}(p^2,-p^2,-p^2)$.

Now consider the homotopy $h:(L'_1 \cup L_2 \cup L_3) \mapsto (L'_1 \cup L'_2 \cup L'_3)$
with support on $L_2 \cup L_3$, for which $L'_1$ is fixed, and the components $L_2, L_3$ are homotopic inside $U$ without self-intersections.
It is convenient to imagine the inverse homotopy $h^{-1}$ as the “pulling apart” of $L'_2 \cup L'_3 $ components located along $\l$ very close to each other into $p$-fold windings
solid tori $V_2 \subset U$, $V_3 \subset U$. Pulling apart occurs without self-intersection, but 
with mutual intersection.

The algebraic number of intersections of components under the homotopy $h^{-1}$ is equal to zero. Let's slow down the crossing
component under homotopy by extending thin whiskers from the $L'_3$ component to the $L'_2$ component. 
Each whisker $\psi_i$ is equipped with a sign $\pm$, which is defined as the coefficient of whisker engagement with $L'_2$.
We denote the whisker family by $\Psi = \{\psi_i\}$.

Family conversion
$\Psi$ into a family $\Psi'$ with the same number of whiskers will be called elementary if the following transformations occur.

A: $\psi_i$ is converted to $\psi'_i$ on a switch with one of the components $L'_3$ or $L'_2$
without twisting and without changing the linking number of the whisker with $L'_2, L'_3$, and the rest of the whiskers do not change. (Note that the linking coefficient of the whisker with $L'_2$, $L'_3$ is defined, since any two whiskers with the same bases and identical points of engagement with $L'_2$ are homotopic.)

B. Cancelation of a close parallel pair of whiskers with opposite signs by isotopy $L_3$.

C. The intersection of  $\psi_i$ with another whisker $\psi_j$ with a sign. (Note that the sign of the intersection of the two whiskers is correctly defined)

D. Self-intersection of the whisker with itself, i.e. twisting 
with a sign. (Note that the sign of the whisker twist (self-intersection) is correctly defined.)

Let's study how $\beta$,$\gamma$-invariants change under transformations A, B, C, D. During transformations
A,B invariants $\gamma$,$\beta$ do not change. Under the transformation B,C, the component $L_2$ remains fixed; any such elementary transformation is decomposed into a composition of homotopies of the component $L_3$.
Using the jump formula for the $\gamma$-invariant, we conclude that
the value of the invariants depends only on the product of the whisker signs (in the case of D, the square of the whisker sign is equal to $+1$) and the sign of the intersection (in the case of $D$, the sign of the self-intersection).   In this case, the change in $\gamma$ is taken with the coefficient $p^2$, and the contribution of the jumps $\gamma$ and $\beta$ to the value of the $\tilde{M}$-invariant is reduced.

Thus, it has been proven that the transformation $\Psi \mapsto \Psi'$ induces a transformation of the link that preserves the value of the $\tilde{M}$-invariant. By a family of elementary formations, any two families of whiskers with the same bases and identical coefficients of engagement are translated into another.  The jump of the $\tilde{M}$-invariant does not depend on exactly how the $\Delta$-motions connecting the indicated links were chosen. Lemma \ref{l23} is proved. \qed

\subsubsection*{Proof of the Lemma \ref{l24}}
Consider the special link $\L_{Hopf}^{op(1)}(p,-p,-p)$. The calculation of $\tilde{M}(\L_{Hopf}^{op(1)}(p,-p,-p))$ is similar to Lemma \ref{Hopf3}. We will untie this link by self-intersections of the components: $\L_{Hopf}^{op(1)}(p,-p,-p)=\L \mapsto \L_1$, and we will fix the sublink $(L_1,L_2) \subset \L$, which we will join by a homotopy with the trivial link; $p$ intersections will be required.

At each intersection during smoothing, a 2-component link $(L_1\sharp L_2\cup L_3)$ arises, which has two disjoint unknotted components, and therefore each time 
the jump $c_1$ is zero.

As a result of untying, an engagement will be obtained, composed of two parallel unknotted (oppositely oriented and unengaged)
a component, when untied, only Hopfian engagements with opposite turns arise.
Therefore we get $c_1(\L_{Hopf}^{op(1)}(p,-p,-p))=0$. By calculation of the invariant  $\gamma$, let us remark that the only non-trivial term is $-(-2k)c_1(L_2 \cup L_3) = \frac{p^4-p^2}{3}$.  
As the result we get: $\tilde{M}(\L_{Hopf}^{op(1)}(p,-p,-p))= \frac{(p^4-p^2)(-p^3)}{3} + \frac{p^4(p^3-p)}{6} = \frac{-p^7+p^5}{6}$.

 By Lemma \ref{Hopf6} the values $\tilde{M}(\L_{Hopf}(p,p,p))$ and $\tilde{M}(\L_{Hopf}^{op(1)}(p,-p,-p))$ are related by the formula:
 $$\tilde{M}(\L_{Hopf}^{op(1)}(p,-p,-p)) - \tilde{M}(\L_{Hopf}(p,p,p)) = - \frac{2p^7}{3} - \frac{p^5}{3}. $$ 
As the result we have:
$\tilde{M}(\L_{Hopf}(p,p,p))=\frac{p^7}{2}+\frac{p^5}{2}$.
\qed

 \subsubsection*{Normalizations of links}

Let us define a normalization  
$\L_{norm}=((2,3)L_1, (3,1)L_2, (1,2)L_3)$, 
for which all pairwise linking coefficients are equal to  $k$, $k=(1,2)(2,3)(3,1)$.
By the construction the standard moidel of $\L_{norm}$ is a special link. 
 Put  $k_{norm}=k^2$.
 The normalized link $\L_{norm}$ by the turns $(\lambda_1,\lambda_2,\lambda_3)$
 is transformed as follows:
$((2,3)L_1, (3,1)L_2, (1,2)L_3) \mapsto (\lambda_1(2,3)L_1, \lambda_2(3,1)L_2, \lambda_3(1,2)L_3)$.  As a result of the turns
 we have
$k_{norm} \mapsto k_{norm}\lambda_1^2\lambda_2^2\lambda_3^2$. We assume that $k$ admits the identity transformation by  turns.

  


Let us consider a good link $\L$, for which the pairwise linking coefficients are equal to  squares and consider ita normalization $\L_{norm}$.  
Let us define the secondary operation $\L_{norm} \mapsto \L_{NORM}$ of taking cables for each component of the link with the multiplicity  $\sqrt{k}=\mu_1^2\mu_2^2\mu_3^2$.

By the cable construction $\L \mapsto (\lambda_1,\lambda_2,\lambda_3)\L$ of a good link its secondary normalization  $\L_{NORM}$ is transformed by the secondary cable construction with the vector
of turns $(\lambda_1\sqrt{\lambda_1\lambda_2\lambda_3},\lambda_2\sqrt{\lambda_1\lambda_2\lambda_3},\lambda_3\sqrt{\lambda_1\lambda_2\lambda_3})$.  Pairwise linking coefficients of a link $\L_{NORM}$ equal to $k_{norm}=k^2$. 

The cable construction defined above can be described by the following dyagram:
\begin{eqnarray}\label{diagK}
	\L  \mapsto  \L_{norm} \mapsto \L_{NORM}.  
\end{eqnarray}

Consider the following invariants: $k^2\tilde{M}^{av}(\L)$,  $k^{-2}k^{-3}_{norm}\tilde{M}^{av}(\L_{NORM})=k^{-8}\tilde{M}^{av}(\L_{NORM})$. 
The invariants admits asymptotic property (\ref{ass}). The invariants
have a same jump under $\Delta$-motions of $\L$, inducing $\Delta$-motions
of $\L_{NORM}$.

 The invariants $k^{-2}k^{-3}_{norm}\tilde{M}^{av}(\L_{NORM})$, $\tilde{M}(\L)^{av}$ are related as follows:
\begin{eqnarray}\label{Rnorm}
k^{-2}k^{-3}_{norm}\tilde{M}^{av}(\L_{NORM}) = \tilde{M}^{av}(\L) + R((2,3),(3,1)(1,2)), 
\end{eqnarray}
where $R$ is a polynomial, which is explicitely defined below.

\begin{lemma}\label{new}
	For a good link $\L$ the polynomial  $R$ in ($\ref{Rnorm}$) (comp. with (\ref{Poly})) is given by 
	\begin{eqnarray}\label{R}
		\begin{array}{c}  
R((2,3),(3,1),(1,2))=\frac{1}{24}(1,2)^3(2,3)^3[(3,1)^3-(3,1)] +\\ 
\frac{1}{24}(1,2)^3(2,3)^3[(3,1)^3-(3,1)]+\frac{1}{24}(1,2)^3(2,3)^3[(3,1)^3-(3,1)].\\
\end{array}
\end{eqnarray}
\end{lemma}

\subsubsection*{Proof of Lemma \ref{new}}
The left- and right- sides of (\ref{Rnorm}) have a same jump by $\Delta$-moves. 
Let us consider the  standard link $\L_0((3,1),(2,3),(1,2))$ with a prescribed collection of pairwise linking numbers. Take $\tilde{M}$ in the right-hand side of the formula.  Let us prove the statement (A): the left-hand side of the formula equals zero.

Take the link $\L_{NORM}$ and consider the sequence of $\Delta$-moves, which joins the link
$\L_{NORM}$ with $\L_0(k^2,k^2,k^2)$. The sequence is decomposed in a collection of elementary $\Delta$-mooves, which permutes  wriskers of component $[L_1,L_2][L_1,L_3]\mapsto [L_1,L_3][L_1,L_2]$. An single commutator contatins $2$ $\delta$-moves, each move has a order $1$ jump of invariants $k^6\gamma$, $k^8\beta$ with respect to linking numvber of wriskers, which depends of its position. The pair of moves has a constant jump, which equals to to 
$\pm k^8$ with opposite signes. This proves the statment (A).
 
Let us consider  the link $\L^{op(3+1)}(-(2,3),(3,1),-(1,2))$ and a sequence of $\Delta$-moves from this link to $L_0(-(2,3),(3,1),-(1,2))$. This sequense induces the sequence of $\Delta$-moves of the corresponding normalized links. Let us prove that the jump of the indused sequence equals to the the first term in the formula (\ref{R}). The jump of the invariant $\tilde{M}$ equals to $\frac{1}{6}(1,2)^2(2,3)^2[(3,1)^3-(3,1)]$. A jump of the induced sequence  of moves of secondary link $\L_{NORM}$ contains a collection of $\sqrt{k}^3$ copies of $\Delta$-moves (involdev 2 components) of primery link $\L_{norm}$. The full collection hase a jump $\tilde{M}^{ev}(\L_{NORM}) \equiv k^6$, each moves is tacken with the coefficient $k^{-6}$. We collect only moves of $\L_{NORM}$, $L_{norm}$ with
components are linked with no opposition of the linking number  independently (a $\frac{1}{4}$  of the full number). A collection of moves of the link $\L_{norm}$ contains $(1,2)(2,3)k^2$ copies of the coresponding moves of
$\L_0$, each jump of $\tilde{M}^{ev}(\L_{norm})$ is tacken with the coefficient $\equiv k^{-2}$.
	
The full collection of jumps from $\L^{op(3+1)}$ into $\L_0((-(2,3),(3,1),-(1,2))$ gives the jump of the invariant in the right-hand side of the formula, which is given by the terms of (\ref{R}). Lemma 	\ref{new} is proved. \qed

\subsubsection*{Internal symmetry of special links}

Let us define an equivalent classes of oriented ordered links with respect to the inversions of orientations of components and transpositions of orders. This equivalent classes are called the internal symmetry classes.

Consider an arbitrary link
$\L=(L_1,L_2,L_3)$. Assume that $\L$ is a right link: $k=(1,2)(2,3)(3,1)>0$. Let us assume that all paiwise linking numbers are positive. In its $\Delta$-equivalen class this link is represented by the 
special link $\L_0((3,1),(1,2),(2,3))$ (see Fig.3). 
 In a general case with the condition $k>0$ we have additional 3 cases when two of 3 linking numbers $(i+1)<0,(i+2)<0$ are negative and one linking number $(i,i+1)>0$ is positive. In this case we take the standard model of $\L_0$ with prescribed collection of the linking numbers, such that $L_{i+2}$ coinsids with the  $i+2$-th component  of the special link $L_0$ and is opposite to the $i$-th and $i+1$-th components
of $\L_0$. This modification of the special link is denoted by $\L_0^{tw}$.

In the case $k<0$ we have the special link with all linking numbers are negative and we have additional 3 cases when two of 3 linking numbers $(i+1)<0,(i+2)>0$ are positive and one linking number $(i,i+1)<0$ is negative. In this case we take the standard model of $\L_0$ with prescribed collection of the linking numbers, such that $L_{i+2}$ coinsids with the inversion of the $i+2$-th component  of the special link $L_0$ and coinsids to the $i$-th and $i+1$-th components of $\L_0$.

\begin{definition}
Drefine the invariant $M$ of a good link $\L$ in (\ref{m}) by the formula:
\begin{eqnarray}\label{MM}
 M(\L) = k^{-2}k^{-3}_{norm}\tilde{M}^{av}(\L_{NORM}),
\end{eqnarray}
where $\L_{NORM}$ is a secondary normalization of $\L$. The right-hand side of (\ref{Rnorm})
is defined without an assumption that $\L$ is good.
\end{definition}

\subsection*{Proof of Theorem \ref{th1}}

Let us start with the statement 1 of the theorem. This is a corollary of the following lemma.

\begin{lemma}\label{l16}
	
1.	The invariant $M(\L)$ for a good link $\L$ satisfy the asymptotic property $(\ref{ass})$.

2.  The invariant $M(\L)$ for an arbitrary link $\L$ satisfy the asymptotic property $(\ref{ass})$.
\end{lemma}

\subsubsection*{Proof of Lemma \ref{l16}}

Statement 1 is a corollary of the formula \ref{eql24} and Lemma \ref{l23}. We apply the formula for $k(\L_{Hopf}(+1,+1,+1))$ and gets the required asymptotic property for $p=k^2$. 

To prove Statement 2 let us show that the asymptotic property holds for an arbitrary link $\L$, and not only under the assumption that the link is good.

Let us define the polynomial $Q(\lambda_1L_1 \cup \lambda_2L_2 \cup \lambda_3L_3)$ depending on the pairwise linking coefficients $\L$ and the vector $(\lambda_1,\lambda_2,\lambda_3)$ according to the formula:
$$\tilde{M}^{av}(\lambda_1L_1 \cup \lambda_2L_2 \cup \lambda_3L_3)-(\lambda_1\lambda_2\lambda_3)^4
\tilde{M}^{av}(L_1 \cup L_2 \cup L_3).$$
For a good link $\L$ we obtain: $Q(\L)=0$. On the other hand, for an arbitrary $\L$, $Q$ depends only on its pairwise
linkage coefficients and is a polynomial.

Let us arrange the monomials of the polynomial $Q$ in lexicographic order, let $\xi(\alpha,\beta,\gamma)(2,3)^{\alpha}(3,1)^{\beta}( 1,2)^ {\gamma}$ be a principal monomial. Let us choose a good link $\L_1=(L_1 \cup L_2 \cup L_3)$ for which $(2,3)^{\alpha}(3,1)^{\beta}(1,2)^{ \gamma } \neq 0$. For example, we can define $\L=Hopf(+1,+1,+1)$. Let's take a look
 $\L=Hopf(+1,+1,+1)$. Concider the vector
$(\lambda_1=\mu_2^2\mu_3^2, \lambda_2=\mu_3^2\mu_1^2, \lambda_1=\mu_1^2\mu_2^2 )$, where
$\mu_3 >> \mu_2 >> \mu_1 > 0$. Then the principal monomial in  $Q(\L,\lambda_1,\lambda_2,\lambda_3)$ is greather (by its absolute value) then
the sum of all last nonomial and we have:
 $Q(\lambda_1L_1,\lambda_2L_2,\lambda_3L_3) \neq 0$. On the other hand, $(\lambda_1L_1,\lambda_2L_2,\lambda_3L_3)$ is a good link, if $\L$ is good link.
 This is possible only if
 $Q \equiv 0$. The asymptotic property is proved for an arbitrary link  $\L$. Lemma \ref{l16}
 and Statment 1 of Theorem \ref{th1} is proved. \qed

\subsubsection*{Proof of Theorem \ref{th1}, Statement 2}
This follows from Lemma \ref{Hopf3}. \qed

\subsubsection*{Proof of Theorem \ref{th1}, Statement 3} 
The both invariants has same jumps by $\Delta$-moves, see \cite{A},\cite{A2} for jumps of the integral invariant. For good links this is a corollary of asymptotic property for the link $Hopf(+1,+1,+1)$. For an arbitrary link this follows from the fact that the considered invariants are finite-type invariants (comp. with Statment 2, Lemma \ref{l16}). \qed   

\subsubsection*{Proof of Theorem \ref{th1}, Statement 4,6} 
 Proofs are clear. \qed

 \subsubsection*{Proof of Theorem \ref{th1}, Statement 5} 
 The both sides of (\ref{Rnorm}) are not changed with respect to the translations. \qed
 

\section*{Calculations}

Consider the Hopf link $\L_{Hopf}(+1,+1,+1) = (L_1 \cup L_2 \cup L_3)$, presented in Fig. 6.
\[ \]
Fig.6
\[  \]
For this link we have:
$\gamma=1$; $(1,2)=(2,3)=(3,1)=+1$; $\beta_{2,3}=\beta_{3,1}=\beta_{1,2}=0$; $\tilde{M}=-1$,
$P=+\frac{3}{4}$, $\tilde{M}^{av}=\tilde{M}+P=-\frac{1}{4}$, $R=0$, $M=\tilde{M}^{av}=-\frac{1}{4}$.
\[  \]

Consider the Hopf link $\L_{Hopf}(-1,+1,+1)$,  presented in Fig. 7.
\[ \]
Fig.7
\[  \]
For this link we have:
$\gamma=0$; $(1,2)=(3,1)=+1$, $(2,3)=-1$; $\beta_{2,3}=\beta_{3,1}=\beta_{1,2}=0$; $\tilde{M}=0$,
$P=\frac{1}{4}$, $\tilde{M}^{av}=\frac{1}{4}$, $R=0$, $M=\tilde{M}^{av}=\frac{1}{4}$.
\[  \]
 
Consider the cable
 $(2L_1 \cup L_2 \cup L_3)$ assotiated with the link $\L_{Hopf}(+1,+1,+1)$, see Fig. 8.
\[  \] 
Fig.8
\[  \]
For this link we have:
$\gamma=6$; $(1,2)=(3,1)=+2$, $(2,3)=+1$; $\beta_{2,3}=\beta_{3,1}=\beta_{1,2}=0$; 
$P_1=16$, $\tilde{M}=-24$, $\tilde{M}^{av}=-8$, $R=4$, $M=\tilde{M}^{av}+R=-4$.

Consider the cable
 $(4L_1 \cup L_2 \cup L_3)$ assotiated with the link  $\L_{Hopf}(-1,-1,-1)$, see Fig. 9,
\[  \] 
Fig.9
\[  \]
which is a good left link. Let us calculate $M$ by the formula
 (\ref{MM}). We have: $\L_{norm} = (4L_1 \cup 4L_2 \cup 4L_3)$, $k_{norm}=16=2^4$;
$\L_{NORM} = (8L_1 \cup 8L_2 \cup 8L_3)$, $k_{NORM}=258=2^8$. $\tilde{M}=2^{38}$. 
$M=2^{-32}2^{38}=2^6=\frac{4^4}{4}=64$. Let us calculate $M$ by the formula (\ref{m}).
We have (for the simplest cabble position): $\gamma = +50$, $\beta_{1,2}=\beta_{2,3}=\beta_{3,1}=0$, $\tilde{M}=-800$, $P=416$, $M^{av}=-384$, $R=320$,
$M=64$.

Consider the cable
 $( 2L_1 \cup L_2 \cup (-1)L_3)$ assotiated with the link $\L_{Hopf}^{op(3)}(-1,-1,+1)$.
\[  \] 
Fig.10.
\[  \]
For this link we have:
$\gamma=6$; $(1,2)=+2$, $(3,1)=-2$, $(2,3)=-1$; $\beta_{1,3}=-1$, $\beta_{1,2}=\beta_{2,3}=0$; 
$P=-4$, $\tilde{M}=-4$. $\tilde{M}^{av}=-8$. This link is isotopic to the link $(-L_1,-2(L_2),L_3)$, which is in the internal symmetry class of the link (Fig.8) with permuted components. We have: $M=-4$.

Consider the link $\L_{Hopf}(+2,+2,+2)$, see Fig.5 (for $p=+2$, $k=+8$). By the formula (\ref{eql24}) we get $\tilde{M}^{av} = \frac{k^2}{4}=16$.  We have $\gamma=31$,  $\tilde{M}=-104$, $P=72$, $\tilde{M}^{av}=-32$, $R=48$, $M=16$.
By the formula (\ref{MM}) we get $M=16$.

\section{Discussion}

We have proven, using only the apparatus of a Conway polynomial in one variable, that there exists an asymptotic invariant $M_3$. This invariant has a finite order in the Vasiliev sense. The main term in the formula (\ref{MInv}) is of order $7$, the polynomial in the formula (\ref{m}) is of degree $7$. This invariant satisfies the asymptotic condition (\ref{ass}) and is oblique. The polynomial $P$ in the formula (\ref{m}) is defined by the formula (\ref{Poly3}), with is the sum of the formulas (\ref{ppoly1}) and (\ref{R}). The term (\ref{ppoly1}) ensures that the invariant is preserved when the leading components of the link are oriented. This polynomial was specified in \cite{A},
while the supposed one (\ref{R}) is new.

The $M_3$ invariant for a magnetic field with closed magnetic lines coincides with the integral invariant defined in \cite{A}. The integral invariant is ergodic, i.e. this invariant is defined as the average integral value of magnetic lines in a compact region (ideal conducting liquid medium); on the boundary surface of the region, the integral curves are directed tangentially to this surface. The invariant of magnetic lines has a density that is represented by a measurable integrable function and does not change under non-compressive transformations of the region that preserve the magnetic field. The invariant corresponds to applied problems that use the theory of the average magnetic field \cite{A5},\cite{A-D}.

In a subsequent publication, the author plans to transfer the results to the case of the asymptotic invariant $M_5$ \cite{A1} and test the asymptotic property of this invariant using the theory of the Conway polynomial in two variables.
  
The author thanks S.A. Melikhov for discussions and pointing out references \cite{Ci},\cite{T}.

\end{document}